\documentclass[review]{elsarticle}

\usepackage{lineno,hyperref}
\modulolinenumbers[5]

\usepackage{mathtools,amsthm,adjustbox}

\begin{document}

\begin{frontmatter}

\title{Approximate conditions admitted by classes of the Lagrangian ${\cal L}=\frac12\left(-u'^2+u^2\right)+\epsilon^iG_i(u, u^\prime, u^{\prime\prime})$}

\author[]{Sameerah Jamal\corref{mycorrespondingauthor}}
\address{School of Mathematics and Centre for Differential Equations,Continuum Mechanics and Applications,
University of the Witwatersrand, Johannesburg, South Africa}
\ead{Sameerah.Jamal@wits.ac.za}
\cortext[mycorrespondingauthor]{Corresponding author}

\author{Nkosingiphile Mnguni}
 \address{School of Mathematics and Centre for Differential Equations,Continuum
Mechanics and Applications, University of the Witwatersrand, Johannesburg,
South Africa}
\ead{nzmnguni8@gmail.com}

\begin{abstract}
We investigate a class of Lagrangians  that admit a  type of perturbed harmonic oscillator  which occupies a special place in the literature surrounding perturbation theory. We establish explicit and generalized geometric conditions for the symmetry determining equations.  The explicit scheme provided can be followed and specialized for any concrete perturbed differential equation possessing the Lagrangian.  A systematic solution of the conditions 
generate  nontrivial approximate symmetries and transformations. Detailed cases are discussed to illustrate the relevance of the conditions, namely (a) $G_1$ as a quadratic polynomial, (b) the Klein-Gordon equation of a particle in the context of Generalized Uncertainty Principle and (c) an orbital equation from an embedded Reissner-Nordstr\"om black hole. 
\end{abstract}

\begin{keyword}
Approximate Symmetries; orbital equation; uncertainty principle.
\MSC[2010] 22E60; 76M60; 35Q75; 34C20
\end{keyword}

\end{frontmatter}


\section{Introduction}
A Lie symmetry group forms a robust tool  in the analysis of differential equations, primarily because it provides invariant functions which may reduce the order of the equation and lead to the determination of analytic solutions. 
Differential equations which possess  a variational principle or Lagrangian, admit specialized Lie symmetries, called Noether symmetries or divergence symmetries, which in addition to the invariant functions,   leave the action invariant. 
Aside from these classical symmetries, there exist approximate symmetries which are devised from  equations, regarded as perturbed equations, that contain some small parameter $\epsilon$. Within the literature, among several computational techniques for  approximate generators, there are two main formalisms, one proposed by Baikov, Gazizov and Ibragimov \cite{bai2} and the second was presented by Fushchich and Shtelen \cite{bai3}. Thereafter, the concept of  approximate Noether symmetries and conservation laws emerged \cite{a1, a3}. Owing to these developments, many important physical differential equations have been studied successfully, see for instance \cite{a2, pek1, bar2}. Note that unlike exact symmetries,  approximate symmetries do not necessarily form a Lie algebra but rather an  ``approximate Lie algebra" \cite{gaz}.

Two decades ago, a method was devised whereby a known symmetry and its corresponding conservation law of a given partial differential equation can be used to construct a Lagrangian for the equation \cite{r1}. However, in the absence of a  Lagrangian,
there has been significant developments on the derivation of approximate conservation laws. For instance in \cite{r2}, a method based on partial Lagrangians  was introduced to
construct approximate conservation laws of approximate Euler-type equations using approximate Noether-type symmetries. 
In \cite{r3}, Zhang  considered  approximate nonlinear self-adjointness for perturbed PDEs and showed how approximate conservation laws, which cannot be obtained by the approximate Noether's theorem, are constructed. Nevertheless, when a Lagrangian is available, Noether's work is not only more elegant, but also highly efficient, and will always be    the preferred method. As an example of the advantages of approximate Noether symmetries, over other existing methods previously mentioned, a recent study, by one of the authors, found a  geometric connection between the Homothetic algebra of an underlying geometry and the approximated Noether symmetries, that is, if the perturbation terms do not modify the Kinetic energy of regular Lagrangians, approximate symmetries exist if and only if the metric that defines the Kinetic energy, admits a nontrivial Homothetic algebra \cite{r4}.  


The purpose of this paper is three-fold. Firstly, in the following work we stipulate the generalized approximate conditions in the case of a class of perturbed Lagrangians, up to third-order, 
\begin{equation}\label{lagg1}{\cal L}(u, u^\prime, u^{\prime\prime},\epsilon)=\frac12\left(-u'^2+u^2\right)+\epsilon^iG_i(u, u^\prime, u^{\prime\prime}).\end{equation} The Lagrangian defined here has the Latin index $i$ that is restricted to the values $1,2\, \textrm{and}\, 3$ and $u$ is a function of $\phi$. The above approximate  class of Lagrangians and its symmetry generators maintain the specified perturbation order of $\epsilon$. To preserve generality we have not made specific assumptions about the $G_i(u, u^\prime, u^{\prime\prime})$. 
Rather we provide an explicit scheme which can be followed and specialized for any concrete differential equation possessing the Lagrangian (\ref{lagg1}), whereby one may extract further information using a given $G_i(u, u^\prime, u^{\prime\prime})$. 
Our next purpose is to use the generalized conditions to find approximate divergence symmetries for several critical cases of interest. 
Thirdly, the latter will be used to establish the associated approximate first integrals by invoking  Noether's theorem.

Before we begin, it is worth mentioning that there are powerful and fully automated software routines to obtain symmetries that are not approximate, commonly referred to as exact  symmetries, for example \cite{champ,bau,sym}. Eliminating all the perturbed terms in the
 Lagrangian (\ref{lagg1}), leads to the derivation of the oscillation equation
\begin{equation}\label{o1}u^{\prime\prime}+u=0.\end{equation}
It is easily seen that this unperturbed equation  is maximally symmetric and admits the 8-dimensional Lie algebra of exact symmetries $sl(3, R)$ given by
$$\begin{array}{lc}
&X_{0}^1=\partial_\phi, \\
&X_0^2=\sin(2 \phi)\partial_\phi+ \cos(2 \phi)u\partial_u, \\
&X_0^3=\cos(2 \phi)\partial_\phi- \sin(2 \phi)u\partial_u, \\
&X_0^4=\sin( \phi)\partial_u, \\
&X_0^5=\cos( \phi)\partial_u,\\
&X_0^6=u\partial_u,\\
&X_0^7=u\cos(\phi)\partial_\phi-u^2\sin(\phi)\partial_u,\\
&X_0^8=u\sin(\phi)\partial_\phi+u^2\cos(\phi)\partial_u.
\end{array}$$

In a problem with a small perturbation, one may consider the approximate Lie symmetry approach versus the approximate Noether or variational symmetry approach.
 We have chosen here the approximate variational  approach since we shall find, at our disposal, explicit formulae for the approximate symmetry conditions and conservation laws ensured by Noether's theorem (see Sections 3), whose determination is usually sans the use of algebraic and algorithmic software. Comparatively, the approximate Lie method is tedious and involves extra computations. Thus it is immediate and far more efficient to apply the variational approach. In order to illustrate our main results or derived conditions, some examples  are presented in the text. These examples are appropriately chosen, for they are novel in the sense that they have not been subjected to an approximate symmetry investigation. Moreover these examples  involve variational principles in a cosmological and relativistic setting.  
 One case explores the approximate symmetries of an orbital equation that arises when a Reissner-Nordstr\"om black hole is embedded into a Friedman-Robertson-Walker (FRW) space \cite{RNpaper}. To obtain the equation of motion of a planet, it is the norm to rewrite a given metric from the cosmic coordinate system to the Schwarzschild or solar coordinate system and thereafter deduce the geodesic equation. Significantly, such equations have the propensity to show whether or not the orbit of a planet is influenced by the evolution of the universe. As a second case, we investigate the modified Klein-Gordon equation of a  spin-0 particle in the Generalized Uncertainty Principle (GUP) \cite{g1,g2,g3,g4, pal}. In general, as detailed below, the  modified Klein-Gordon equation is a fourth-order partial differential equation, which we reduce and adapt to possess the perturbed Lagrangian (\ref{lagg1}).  In each case, we state the approximate first integrals corresponding to the  approximate Noether symmetries obtained.

\noindent The plan of the paper is as follows. In the next section we briefly review the geometric preliminaries surrounding exact and approximate point symmetries of differential equations, with a focus on generators originating from a variational principle. This section also 
introduces the  notation and conventions assumed.
 The perturbed class of Lagrangians (\ref{lagg1}) are studied in Section 3, where we show that the  approximate symmetry determining equations are generated by a set of generic conditions. 
 In Section 4, we apply the general results of the previous sections to highlight a particular case of $G_1(u, u^\prime, u^{\prime\prime})$ that admits an enlarged ``group" of approximate Noether symmetry generators.
 Section 5 describes the case of  the modified Klein-Gordon equation of a particle in the GUP while Section 6 draws attention to several orbital equations of interest where the generalized conditions are especially useful. Finally, in section 7 we present our conclusions.

\section{Point transformations}

\label{section3}

Our  interest lies in point transformations, and for the convenience
of the reader we insert the necessary theory pertaining to this analysis. The presentation here is for ordinary differential equations, however most of the theory has been generalized to partial differential equations in the references cited. For the sake of brevity, the summation convention is adopted in this text, in which there is summation over all repeated  indices.  First, consider a system of second-order ordinary differential
equations (unperturbed), where $t$ is the independent variable and $x^{i}$ denotes the dependent variables
\begin{equation}
{x}^{\prime\prime i}=\omega^{i}\left(  t,x^{j},{x}^{\prime j}\right). \label{Lie.0}%
\end{equation}

An one-parameter point transformation in the space $\left\{  t,x^{j}\right\}
$,  has the property of mapping solutions of
(\ref{Lie.0}) to itself and satisfies the infinitesimal 
criterion of invariance
\[
X^{\left[  2\right]  }\left(  {x}^{\prime\prime i}-\omega^{i}\right)  =0
\quad\mathrm{mod} \quad{x}^{\prime\prime i}-\omega^{i}=0,
\]
where  $X$ is defined as
\begin{equation}
X=\frac{\partial\bar{t}}{\partial\varepsilon}\Bigg|_{\varepsilon=0}%
\partial_{t}+\frac{\partial\bar{x}^{i}}{\partial\varepsilon}%
\Bigg|_{\varepsilon=0}\partial_{i} \label{Lie.01}%
\end{equation}
with $X^{\left[  2\right]  }$ as the second prolongation of $X$ in the jet
space of variables \cite{StephaniB}. $X$ is the generator of the point transformation  called a\ Lie symmetry for the system of differential equations. 

On the other hand, if the  system (\ref{Lie.0}) follows from the
variation of the action integral $$S=\int\mathcal{L}dt,$$ then the Noether's
theorem  \cite{Noetherpaper} states that when a (finite) group of transformations leaves the
action  invariant, i.e. 
\begin{equation}
\label{Noethcon}S\left(  t,x^{j},...\right)  = S\left(  \bar{t}\left(
t,x^{j},\varepsilon\right)  ,\bar{x}^{j}\left(  t,x^{j},\varepsilon\right)
,...\right)  ,
\end{equation}
then a conserved quantity exists. For an unperturbed Lagrangian  up-to first-order in
derivatives (as is the case in this paper), condition (\ref{Noethcon}) yields
\begin{equation}
X^{\left[  1\right]  }\mathcal{L}+\mathcal{L}\frac{d}{dt}\left(
\frac{\partial\bar{t}}{\partial\varepsilon}\Bigg|_{\varepsilon=0}\right)
={f^\prime} \label{Lie.03}%
\end{equation}
and the corresponding first integral is given by
\begin{equation}\label{fir}
I= \left(  {x}^{\prime j}\frac{\partial\mathcal{L}}{\partial{x}^{\prime j}%
}-\mathcal{L}\right)  \frac{\partial\bar{t}}{\partial\varepsilon
}\Bigg|_{\varepsilon=0} -\frac{\partial\mathcal{L}}{\partial{x}^{\prime j}}
\frac{\partial\bar{x}^{j}}{\partial\varepsilon}\Bigg|_{\varepsilon=0} +f.
\end{equation}
In this scenario $X$ would be called a Noether symmetry which is also a Lie symmetry; however the
inverse of this result can be false.
We continue with  the 
 review of the techniques of finding approximate variational symmetries.
For a $k$-th order perturbed system of ordinary differential equations
\begin{equation}
E=E_0+\epsilon E_1+\epsilon^2 E_2+\ldots+\epsilon^k E_k+O(\epsilon ^{k+1}),
\end{equation}
corresponding to a  Lagrangian, which is perturbed in $ \epsilon $,
\begin{equation}
{\cal L}(t,x,x^{\prime j},\epsilon)={\cal L}_0(t,x,x^{\prime j})+\epsilon {\cal L}_1(t,x,x^{\prime j})+\ldots+O(\epsilon ^{k+1}),
\end{equation}
the functional $ \int {\cal L}dt $ is invariant under the one-parameter group of transformations with approximate Lie symmetry generator
\begin{equation}
X=X_0+\epsilon X_1+\ldots+\epsilon^k X_k,
\end{equation} up to gauge
\begin{equation}
A=A_0+\epsilon A_1+\ldots+\epsilon^k A_k,
\end{equation}
if
\begin{equation} \label{apcon}
X{\cal L}+\left(D_t \frac{\partial\bar{t}}{\partial\varepsilon}\Bigg|_{\varepsilon=0}\right){\cal L}=D_tA,
\end{equation}
where $ D_t $ is the total derivative operator.
In this notation, 
$X_0$ is the exact symmetry generator originating from the unperturbed Lagrangian and $X_1$ the first-order approximate symmetry generator.
A perturbed equation always admits the trivial approximate symmetry generator
$\epsilon X_0$. Also, if $X = X_0 + \epsilon X_1$ exists with $X_0\ne0$ and $X_1\ne kX_0$ (k an arbitrary constant), then it is called a nontrivial symmetry \cite{ma}. These considerations can be generalized to higher-order approximate symmetry generators.
An analogous formula for the first-order approximate first integrals can be obtained from Eq. (\ref{fir}) bearing in mind that for first-order, the approximate first integrals are defined by $I = I_0 + \epsilon I_1$, where $I_0$ is the exact first integral and $I_1$ is the first-order approximate part.

\section{Third-Order Geometric Conditions}

First we obtain the approximate Noether symmetry conditions for the class of Lagrangians (\ref{lagg1}) by applying the  approximate symmetry theory. Then we shall study the approximate Noether symmetries of (\ref{lagg1}) pertaining to several important problems in the literature. The determination of approximate Noether point symmetries of the Lagrangian (\ref{lagg1}) consists of two steps: (a) the derivation of the conditions which provide the symmetry determining equations, and (b) the solution of these determining equations. The first step is precise, however the symmetry conditions which arise can be quite involved. The key point is to express conditions for generic forms of $G_i(u, u^\prime, u^{\prime\prime})$.

As mentioned above, Noether symmetries are just a specialization of Lie symmetries, and thus the $sl(3, R)$ algebra given above contains  the Noether point symmetry generators. The latter  comprises of a 5-dimensional Lie algebra $X_0^{1-5}$ with
the corresponding gauge term (the $c_j$ are constants)
$$A_0={u}^{2}\cos \left( 2\,\phi \right) c_3+{u}^{2}\sin \left( 2\,\phi
 \right) c_2+\sin \left( \phi \right) c_5\,u-\cos \left( \phi
 \right) c_4\,u+c_6.$$
The Noether first integrals  corresponding to each $X_0^h, h=1,\ldots,5$ are
$$\begin{array}{lc}
&I_0^1=\frac{1}{2}\left(u^2+u'^2\right), \\
&I_0^2=\frac{1}{2}\left(u'^2-u^2\right)\sin(2 \phi)- uu'\cos(2 \phi), \\
&I_0^3=\frac{1}{2}\left(u'^2-u^2\right)\cos(2 \phi)+ uu'\sin(2 \phi),\\
&I_0^4=-u'\sin( \phi)+ u\cos( \phi), \\
&I_0^5=-u'\cos( \phi)- u\sin( \phi).
\end{array}$$

If we  include a perturbation up to first-order in $\epsilon$, that is, the Lagrangian (\ref{lagg1})  omits the terms in $G_2(u, u^\prime, u^{\prime\prime})$ and $G_3(u, u^\prime, u^{\prime\prime})$,  the determination of approximate symmetries takes a particular form. That is, for each term of the Noether condition (\ref{apcon}) for the Lagrangian (\ref{lagg1}) we have the geometric condition
\begin{equation}\begin{array}{lc} \label{con1}
&X{\cal L}=\left(\eta_{1,\phi}+u^{\prime}\eta_{1,u}-u^{\prime}\xi_{1,\phi}-\left(u^{\prime}\right)^2\xi_{1,u}\right)\left(-u^{\prime}\right)+\eta_1u+ \\
&\left(-2\,\sin \left( \phi \right) \cos \left( \phi \right) c_3\,u_{{}}+2\,
c_2\,u_{{}} \left( \cos \left( \phi \right)  \right) ^{2}+c_4\sin \left( \phi \right) +c_5\,\cos \left( \phi \right) -c_2u_{{}}
\right)G_{1,u}\\
&+\left(-4\,c_2\,u_{{}}\cos \left( \phi \right) \sin \left( \phi \right) +2\,
\sin \left( \phi \right) \cos \left( \phi \right) c_3\,u^{\prime}-2\,
 \left( \cos \left( \phi \right)  \right) ^{2}c_2\,u^{\prime}\right)G_{1,u'}\\
 &\left(-4\,
 \left( \cos \left( \phi \right)  \right) ^{2}c_3\,u_{{}}-c_5\sin \left( \phi \right) +c_4\,\cos \left( \phi \right) +c_2u^{\prime}+2\,c_3\,u_{{}}
\right)G_{1,u'}\\
&+\left(8\,\sin \left( \phi \right) \cos \left( \phi \right) c_3\,u_{{}}+6\,
\sin \left( \phi \right) \cos \left( \phi \right) c_3\,u^{\prime\prime}-8c_2u_{{}} \left( \cos \left( \phi \right)  \right) ^{2}\right)G_{1,u''}\\
&-\left(6
 \left( \cos \left( \phi \right)  \right) ^{2}c_2\,u^{\prime\prime}-c_4\sin \left( \phi \right) -c_5\,\cos \left( \phi \right)
 +4c_2u_{{}}+3\,c_2\,u^{\prime\prime}
\right)G_{1,u''},\\
\end{array}
\end{equation}
\begin{equation}\label{x1}
\left(D_\phi \frac{\partial\bar{\phi}}{\partial\varepsilon}\Bigg|_{\varepsilon=0}\right){\cal L}=\left(2c_2\cos(2\phi)-2c_3\sin(2\phi)\right)G_1+\left(\xi_{1,\phi}+u^{\prime}\xi_{1,u}\right)\left(-\frac{\left(u^{\prime}\right)^2}{2}+\frac{u^2}{2}\right),
\end{equation}
\begin{equation}\label{g1}
D_\phi A =A_{1,\phi}+u^{\prime}A_{1,u}.
\end{equation}

On the other hand, if the perturbation is up to second-order in $\epsilon$, the Lagrangian (\ref{lagg1})  omits $G_3(u, u^\prime, u^{\prime\prime})$, and in this case the second condition  is:
\begin{equation}\begin{array}{lc}\label{con2}
X{\cal L}=+\left(\eta_{2,\phi}+u^{\prime}\eta_{2,u}-u^{\prime}\xi_{2,\phi}-\left(u^{\prime}\right)^2\xi_{2,u}\right)\left(-u^{\prime}\right)+\eta_2u +\eta_1G_{1,u}\\
+\left(- {u^{\prime}}^{2}\xi_{1,u}+ u^{\prime}\eta_{1,u}- u^{\prime}\xi_{1,\phi}+\eta_{1,\phi}
\right)G_{1,u'}+\left(-{u^{\prime}}^{3}\eta_{1,uu}-2\,  {u^{\prime}}^{2}\xi_{1,u\phi}\right)G_{1,u''}\\
+\left( \eta_{1,uu} {u^{\prime}}^{2}-3\, \xi_{1,u} u^{\prime}u^{\prime\prime}
+2\, \eta_{1,u\phi} u^{\prime}\right)G_{1,u''}\\
-\left(\xi_{1,\phi\phi}
u^{\prime}+ \eta_{1,u} u^{\prime\prime}-2\, \xi_{1,\phi}u^{\prime\prime}
+\eta_{1,\phi\phi}
\right)G_{1,u''}\\
+\left(-2\,\sin \left( \phi \right) \cos \left( \phi \right) c_3\,u_{{}}+2\,
c_2\,u_{{}} \left( \cos \left( \phi \right)  \right) ^{2}+c_4\sin \left( \phi \right) +c_5\,\cos \left( \phi \right) -c_2u_{{}}
\right)G_{2,u}\\
+\left(-4\,c_2\,u_{{}}\cos \left( \phi \right) \sin \left( \phi \right) +2\,
\sin \left( \phi \right) \cos \left( \phi \right) c_3\,u^{\prime}-2\,
 \left( \cos \left( \phi \right)  \right) ^{2}c_2\,u^{\prime}\right)G_{2,u'}\\
-\left(4\,
 \left( \cos \left( \phi \right)  \right) ^{2}c_3\,u_{{}}-c_5\sin \left( \phi \right) +c_4\,\cos \left( \phi \right) +c_2u^{\prime}+2\,c_3\,u_{{}}
\right)G_{2,u'}\\
+\left(8\,\sin \left( \phi \right) \cos \left( \phi \right) c_3\,u_{{}}+6\,
\sin \left( \phi \right) \cos \left( \phi \right) c_3\,u^{\prime\prime}-8c_2u_{{}} \left( \cos \left( \phi \right)  \right) ^{2}\right)G_{2,u''}\\
-\left(6\,
 \left( \cos \left( \phi \right)  \right) ^{2}c_2\,u^{\prime\prime}-c_4\sin \left( \phi \right) -c_5\,\cos \left( \phi \right) +4c_2u_{{}}+3\,c_2\,u^{\prime\prime}
\right)G_{2,u''}
\end{array}\end{equation}
$$\begin{array}{lc}\label{x2}
\left(D_\phi \frac{\partial\bar{\phi}}{\partial\varepsilon}\Bigg|_{\varepsilon=0}\right){\cal L}=&\left(\xi_{1,\phi}+u^{\prime}\xi_{1,u}\right)G_1+\left(\xi_{2,\phi}+u^{\prime}\xi_{2,u}\right)\left(-\frac{\left(u^{\prime}\right)^2}{2}+\frac{u^2}{2}\right)\\ &+\left(2c_2\cos(2\phi)-2c_3\sin(2\phi)\right)G_2,
\end{array}$$
\begin{equation}\label{g2}
D_\phi A =A_{2,\phi}+u^{\prime}A_{2,u}.
\end{equation}

Last but not least, a third-order perturbation in $\epsilon$ results in the third condition
\begin{equation}\begin{array}{lc}\label{con3}
&X{\cal L}=\left(\eta_{3,\phi}+u^{\prime}\eta_{3,u}-u^{\prime}\xi_{3,\phi}-\left(u^{\prime}\right)^2\xi_{3,u}\right)\left(-u^{\prime}\right)+\eta_3u \\
&+\left(-2\,\sin \left( \phi \right) \cos \left( \phi \right) c_3\,u_{{}}+2\,
c_2\,u_{{}} \left( \cos \left( \phi \right)  \right) ^{2}+c_4\sin \left( \phi \right) +c_5\,\cos \left( \phi \right) -c_2u_{{}}
\right)G_{3,u}\\
&+\left(-4\,c_2\,u_{{}}\cos \left( \phi \right) \sin \left( \phi \right) +2\,
\sin \left( \phi \right) \cos \left( \phi \right) c_3\,u^{\prime}-2\,
 \left( \cos \left( \phi \right)  \right) ^{2}c_2\,u^{\prime}\right)G_{3,u'}\\
 &-4\left(\left( \cos \left( \phi \right)  \right) ^{2}c_3\,u_{{}}-c_5\sin \left( \phi \right) +c_4\,\cos \left( \phi \right) +c_2u^{\prime}+2\,c_3\,u_{{}}
\right)G_{3,u'}\\
&+\left(8\,\sin \left( \phi \right) \cos \left( \phi \right) c_3\,u_{{}}+6\,
\sin \left( \phi \right) \cos \left( \phi \right) c_3\,u^{\prime\prime}-8c_2u_{{}} \left( \cos \left( \phi \right)  \right) ^{2}\right)G_{3,u''}\\
&-6\left(\left( \cos \left( \phi \right)  \right) ^{2}c_2\,u^{\prime\prime}-c_4\sin \left( \phi \right) -c_5\,\cos \left( \phi \right) +4c_2u_{{}}+3\,c_2\,u^{\prime\prime}
\right)G_{3,u''}\\
&+\eta_1G_{2,u}+\left(- {u^{\prime}}^{2}\xi_{1,u}+ u^{\prime}\eta_{1,u}- u^{\prime}\xi_{1,\phi}+\eta_{1,\phi}
\right)G_{2,u'}\\&+\left(-{u^{\prime}}^{3}\eta_{1,uu}-2\,  {u^{\prime}}^{2}\xi_{1,u\phi}\right)G_{2,u''}\\
&+\left( \eta_{1,uu} {u^{\prime}}^{2}-3\, \xi_{1,u} u^{\prime}u^{\prime\prime}
+2\, \eta_{1,u\phi} u^{\prime}\right)G_{2,u''}-\\&\left(\xi_{1,\phi\phi}
u^{\prime}+ \eta_{1,u} u^{\prime\prime}-2\, \xi_{1,\phi}u^{\prime\prime}
+\eta_{1,\phi\phi}
\right)G_{2,u''}\\
&+\eta_2G_{1,u}+\left(- {u^{\prime}}^{2}\xi_{2,u}+ u^{\prime}\eta_{2,u}- u^{\prime}\xi_{2,\phi}+\eta_{2,\phi}
\right)G_{1,u'}\\&+\left(-{u^{\prime}}^{3}\eta_{2,uu}-2\,  {u^{\prime}}^{2}\xi_{2,u\phi}\right)G_{1,u''}\\
&+\left( \eta_{2,uu} {u^{\prime}}^{2}-3\, \xi_{2,u} u^{\prime}u^{\prime\prime}
+2\, \eta_{2,u\phi} u^{\prime}\right)G_{1,u''}\\&-\left(\xi_{2,\phi\phi}
u^{\prime}+ \eta_{2,u} u^{\prime\prime}-2\, \xi_{2,\phi}u^{\prime\prime}
+\eta_{2,\phi\phi}
\right)G_{1,u''},\\
\end{array}
\end{equation}
\begin{equation}\begin{array}{lc}\label{x3}
\left(D_\phi \frac{\partial\bar{\phi}}{\partial\varepsilon}\Bigg|_{\varepsilon=0}\right){\cal L}=&\left(\xi_{2,\phi}+u^{\prime}\xi_{2,u}\right)G_1+\left(2c_2\cos(2\phi)-2c_3\sin(2\phi)\right)G_3+\\&\left(\xi_{1,\phi}+u^{\prime}\xi_{1,u}\right)G_2+
\left(\xi_{3,\phi}+u^{\prime}\xi_{3,u}\right)\left(-\frac{\left(u^{\prime}\right)^2}{2}+\frac{u^2}{2}\right),
\end{array}
\end{equation}
\begin{equation}\label{g3}
D_\phi A =A_{3,\phi}+u^{\prime}A_{3,u}.
\end{equation}

We remark that the above three conditions must be applied sequentially. The separation and solution of the conditions (\ref{con1})-(\ref{g3}) gives the approximate coefficients of the Noether point symmetry vectors. 
To obtain a group classification  involving a generic function, in our case $G_i$, consists of finding the approximate point symmetries of the given Lagrangian with arbitrary $G_i$, and, thereafter to determine all possible and particular cases of $G_i$ for which the symmetry group can be expanded. Naturally there should be a geometrical or physical motivation  in place for considering such specific cases \cite{o}. 
In the following sections we proceed with the applications of conditions (\ref{con1})-(\ref{g3}) in  cases of special interest, that is we deal with some equations admitted by the class of Lagrangians (\ref{lagg1}). Specifically, we study the approximate point symmetries of the modified Klein-Gordon under GUP  and secondly, the approximate point symmetries of an orbital equation arising from an embedded Reissner-Nordstr\"om black hole.
The presentation of results is schematic so not to increase the volume of the paper.

\label{sec:level1}

\section{A Quadratic Polynomial $G_1(u, u^\prime, u^{\prime\prime})$}
If we allow the conditions (\ref{con1})-(\ref{g1}) to act as a selection rule for the functional form of $G_1(u, u^\prime, u^{\prime\prime})$,  we find some surprising  results. In fact the notion of using a symmetry as selection criteria for the free functions  or  recursion operator of a model, can be traced to many articles (for example \cite{prd, prd2, sit} or \cite{olv2, bmm}, respectively). 
The specific case in which $$G_1(u, u^\prime, u^{\prime\prime})=\frac12a_0u^2+a_1u+a_2,\quad$$ we find the exact symmetry algebra 
$X_0^{1-5}$ plus the added approximate generators
$$\begin{array}{lc}
&\bar{X}_{M}^{\epsilon}=X_0^2+\epsilon\bigg[\left(a_0\phi \cos(2\phi)- a_0\sin(2\phi)-2\cos^2(\phi)a_1+a_1\right)\partial_\phi-\\&\left(\frac12a_0u\cos(2\phi)+ua_0\phi\sin(2\phi)\right)\partial_u\bigg],\\
&\bar{X}_{N}^{\epsilon}=X_0^3+\epsilon\left[\left(-\frac12a_0\cos(2\phi)- \phi a_0\sin(2\phi)\right)\partial_\phi-\left(\phi a_0u\cos(2\phi)+a_1\sin(2\phi)\right)\partial_u\right],\\
&\bar{X}_{O}^{\epsilon}=X_0^4+\epsilon\left(\frac12a_0\phi\cos(\phi)-\frac14a_0\sin(\phi)\right)\partial_u,\\
&\bar{X}_{P}^{\epsilon}=X_0^5-\epsilon\left(\frac14a_0\cos(\phi)+\frac12\phi a_0\sin(\phi)\right)\partial_u.
\end{array}$$
Clearly, form the arbitrary constants in $a_m$, the  $a_0$ and $a_1$ should be  nonzero to maximize the number of possible approximate generators.
Correspondingly, the first-order approximate gauge term in this case is
\begin{equation}\begin{array}{lc}
\bar{A}_1=&\frac14\, \left(  \left(  \left( 4\,c_{2}\phi+2\,c_{3} \right) c_6-2\,c_{10} \right) {u_{{}}}^{2}+8\,c_{3}\,c_7\,u_{{}}+4\,c_{3}\,a_{2} \right) \cos \left( 2\phi \right) +\\
&\frac14\, \left(  \left( -4\,c_{3}\,a_{0}\phi+2\,c_{9}
 \right) {u_{{}}}^{2}+8\,c_{2}\,a_{1}\,u_{{}}+4\,c_{2}
\,a_{2} \right) \sin \left( 2\phi \right) +\\
&\frac14\, \left(  \left( 
 \left( 2\,c_{5}\phi-c_{4} \right) a_{0}-4\,c_{13}
 \right) u_{{}}-4\,c_{4}\,a_{1} \right) \cos \left( \phi
 \right)+\\
 &\frac14\, \left(  \left(  \left( 2\,c_{4}\phi+c_{5}
 \right) a_{0}+4\,c_{12} \right) u_{{}}+4\,c_{5}\,c_7 \right) \sin \left( \phi \right) +c_{14}
\end{array}
\end{equation}
The first-order approximate first integrals related to $\bar{X}_{M-P}^{\epsilon}$ are given by
\begin{eqnarray*}
\bar{I}^M=I_0^2+\epsilon\Bigg(-2\,\sin \left( \phi \right) \cos \left( \phi \right) a_1\,u_{{}}-a_0\,\phi \left( \cos \left( \phi \right)  \right) ^{2}{u_{{}}}^{2}+a_0\,\phi \left( \cos \left( \phi \right)  \right) ^{2}{u'}^{2}+\\
1
/2\,a_0\,\phi{u_{{}}}^{2}-1/2\,a_0\,\phi{u'}^{2}-a_0\,\sin \left( \phi \right) \cos \left( \phi \right) {u'}^{2}+u'a_0\,u_{{}} \left( \cos \left( \phi \right)  \right) ^{2}-1/2\,u'a_0\,u_{{}}+\\
2\,u'u_{{}}a_0\,\phi\sin \left( \phi
 \right) \cos \left( \phi \right) -2\,u' \left( \cos \left( \phi
 \right)  \right) ^{2}a_1+u'a_1
\Bigg),\\
\bar{I}^{N}=I_0^3+\epsilon\Bigg(-u'\phi a_0\,u_{{}}+2\,u'a_1\,\sin \left( \phi
 \right) \cos \left( \phi \right) +\\
 \phi a_0\,\sin \left( \phi \right) 
\cos \left( \phi \right) {u_{{}}}^{2}-\phi a_0\,\sin \left( \phi \right) 
\cos \left( \phi \right) {u'}^{2}+2\,u'\phi a_0\,u_{{}}
 \left( \cos \left( \phi \right)  \right) ^{2}-\\
 2\, \left( \cos \left( \phi
 \right)  \right) ^{2}a_1\,u_{{}}-1/2\, \left( \cos \left( \phi
 \right)  \right) ^{2}a_0\,{u_{{}}}^{2}-1/2\,a_0\,
 \left( \cos \left( \phi \right)  \right) ^{2}{u'}^{2}\\
 +1/4\,a_0\,{u'}^{2}+1/4\,a_0\,{u_{{}}}^{2}+a_1\,u_{{}}
\Bigg),\\
\bar{I}^{O}=I_0^4+\epsilon\bigg(-1/2\,u'a_0\,\phi\cos \left( \phi \right) +1/4\,u'a_0\,\sin \left( \phi \right) +1/4\,\cos \left( \phi \right) a_0\,u_{{}
}+\\
\cos \left( \phi \right) a_1-1/2\,a_0\,\phi u_{{}}\sin
 \left( \phi \right) 
\bigg),\\
\bar{I}^{P}=I_0^5+\epsilon\bigg(1/4\,u'a_0\,\cos \left( \phi \right) +1/2\,u'\phi a_0
\,\sin \left( \phi \right) -1/2\,\phi a_0\,u_{{}}\cos \left( \phi
 \right) -\\
 1/4\,\sin \left( \phi \right) a_0\,u_{{}}-\sin \left( \phi
 \right) a_1
\bigg).
\end{eqnarray*}
Next, we apply the conditions to some 
particularly relevant  physical choices of $G_i(u, u^\prime, u^{\prime\prime})$.

\section{The modified Klein-Gordon equation under GUP}

The modified Klein-Gordon equation,  is a fourth-order partial differential equation
\begin{equation}\label{a11}
\Delta\Psi-2\beta h^2\Delta(\Delta\Psi)+V_0\Psi=0
\end{equation}
where $ V_0=\left(\frac{mc}{h}\right)^2 $, $\Delta$ is the Laplace operator 
and the terms $ O\left(\beta^{2}\right) $ have been ignored.\\
The action of the modified Klein-Gordon equation (\ref{a11}) is
$$S=\int dx^4 \sqrt{-g} {\cal L}_A(\Psi, D_\sigma \Psi ),$$
where the Lagrangian $ {\cal L}_A(\Psi, D_\sigma \Psi )$ is given by
\begin{equation}\label{gupl}
{\cal L}_A=\frac{1}{2}\left(  \sqrt{-g} g^{\mu\nu}D_\mu \Psi D_\nu \Psi-  \sqrt{-g}V_0\Psi^2\right).
\end{equation}
Changing variables $\Psi\equiv u$ and reducing Eq. (\ref{a11}), we have the reduced Klein-Gordon equation with $V_0=1, \epsilon=-2\beta h^2$,  that is a fourth-order ordinary differential equation, 
which then possesses the Lagrangian  Eq. (\ref{gupl})  rewritten in the form (\ref{lagg1}),  with $$G_1(u, u^\prime, u^{\prime\prime})=-\frac12 (u^{\prime\prime})^2\quad \textrm{and}\quad G_2(u, u^\prime, u^{\prime\prime})=G_3(u, u^\prime, u^{\prime\prime})=0.$$ 

 After the application of the conditions (\ref{con1})-(\ref{g1}) we find a system of five equations after separation of monomials. The resultant symmetries are $X_0^{1-5}$ and thus the modified Klein-Gordon equation  under GUP contains no first-order nontrivial approximate symmetries.  

On a side note, in this case a trivial approximate symmetry $X^\epsilon=\epsilon \partial_t$  would give the approximate obvious first-order first integral $I=\frac12\epsilon\left(u'^2+u^2\right)$.

\section{The Radial Orbital equation }

The orbital equation  or
 motion equation of a planet, from an embedded Reissner-Nordstr\"om black hole is given by
\begin{equation}\label{a2}
u^{''}+u=\frac{M}{L^2}-\frac{Q^2u}{L^2}+3Mu^2-2Q^2u^3-\frac{H^2}{L^2u^3},
\end{equation}
where $\displaystyle{u\equiv \frac{1}{r}}$, the prime denotes differentiation with respect to $\phi$ and $L$ is the angular momentum of the planet. The terms $3Mu^2$ and $\frac{H^2}{L^2u^3}$ come from the
general relativity and cosmic expansion effect, respectively. Furthermore, the term $\frac{Q^2u}{L^2}$ and $-2Q^2u^2$ are related to  charge. 
The ratio between the term in $H$ and $\frac{M}{L^2}$ is $8 \times 10^{-34} \, \textrm{for Mecury}$ and
$3.6 \times 10^{-28} \, \textrm{for Neptune}.$ If we choose $$\epsilon= 2M\quad\textrm{and}\quad \kappa\epsilon^2=Q^2 \quad\textrm{and}\quad \rho\epsilon^2=H^2,$$
the Lagrangian corresponding to Eq.(\ref{a2}) is given by the general Lagrangian  (\ref{lagg1}) with
\begin{equation}\label{lag2}
G_1(u, u^\prime, u^{\prime\prime})=\left(-\frac{u}{2L^2}-\frac{u^3}{2}\right)\quad \textrm{and}\quad G_2(u, u^\prime, u^{\prime\prime})=\left(\frac{\kappa u^2}{2L^2}+\frac{2\kappa u^4}{4}-\frac{\rho}{2L^2}u^{-2}\right).
\end{equation}

The first step is to retain the term in $G_1(u, u^\prime, u^{\prime\prime})$ from Eq.(\ref{a2}). Consequently, the conditions   (\ref{con1})-(\ref{g1}) provide a system of four equations that solve to give  the first-order approximate Noether symmetry generators
given by
$$\begin{array}{lc}
&X^{\epsilon}_{1}=X_0^4+\epsilon\left(2\sin(\phi)\partial_\phi+u\cos(\phi)\partial_u\right),\\
&X^{\epsilon}_{2}=X_0^5-\epsilon\left(2\cos(\phi)\partial_\phi-u\sin(\phi)\partial_u\right).
\end{array}$$
Correspondingly, the first-order approximate gauge term in this case is
$$\begin{array}{lc}
A_1=&\frac{1}{2L^2}\big(-\cos\left(2\,\phi\right)c_7\,{u}^{2}{L}^{2}+\sin\left(2\,\phi\right)c_6\,{u}^{2}{L}^{2}\\
&+ \left(\left(-c_4\,{u}^{2}-2\,c_9\,u\right){L}^{2}+c_4\right)\cos\left(\phi\right)\big)\\
&+\frac{1}{2L^2}\left(\left(\left(c_5\,{u}^{2}+2\,c_{10}\,u\right) {L}^{2}-c_5\right)\sin\left(\phi\right)+2\,c_{11}\,{L}^{2}\right)
\end{array}$$
The first-order approximate first integrals related to $X_{1-2}^{\epsilon}$ are given by
$$\begin{array}{lc}
I_1^1=&I_0^4+\epsilon\left(\sin(\phi)u'^2-\cos(\phi)uu'+\frac{1}{2}\frac{(L^2u^2+1)}{L^2}\sin(\phi)\right),\\
I_1^2=&I_0^5-\epsilon\left(u'^2\cos(\phi)+uu'\sin(\phi)+\frac12\frac{L^2u^2+1}{L^2}\cos(\phi)\right).
\end{array}$$

In the second approximation, we retain the quadratic 
 $\epsilon$ terms, that is the  $G_1(u, u^\prime, u^{\prime\prime})$ and $G_2(u, u^\prime, u^{\prime\prime})$ defined for  Eq.(\ref{a2}). We proceed with the consideration of the conditions (\ref{con2})-(\ref{g2}) and observe that Eq.(\ref{a2}) possesses no nontrivial second-order approximate symmetry generators, but the first-order approximate symmetry generators are preserved.

\section{Discussion and Conclusion}
In this work we studied the approximate Noetherian point symmetries  and first integrals of the class of  differential equations which follow from a
 Lagrangian perturbed  up-to third-order in $\epsilon$.  We presented new examples where the application of our conditions can be seen. The knowledge of approximate symmetries was used to find the approximate first integrals of the corresponding approximate  equations. 
 We believe that this work can be very useful in the study of various differential  problems. Indeed numerous equations  originate from the generalized Lagrangian (\ref{lagg1}), such as the orbital equations of perturbed spaces. Conditions (\ref{con1})-(\ref{g3}),  applied to the problems studied in  \cite{bar,sch,rn}, immediately gives the results on approximate symmetries (Table \ref{tab:Table_I}), obtained in these works, for orbital equations. 
 \begin{table}[ht!]
\caption{Approximate generators derived using conditions  (\ref{con1})-(\ref{g3}) for several interesting  equations}%
\label{tab:Table_I}%
\begin{adjustbox}{max width=\linewidth}
\begin{tabular}
[c]{c|c|c|c}\hline\hline
\textrm{Model  or} &$G_1(u, u^\prime, u^{\prime\prime})$  &
$G_2(u, u^\prime, u^{\prime\prime})$ &$G_3(u, u^\prime, u^{\prime\prime})$\\
Orbital equation\footnotemark & \textrm{\& Approx. Symmetry}  &
 \textrm{\& Approx. Symmetry}   & \textrm{\& Approx. Symmetry} \\\hline
\textrm{Schwarzschild} &$\left(-\frac{u}{2L^2}-\frac{u^3}{2}\right)$ & $0$ &
0\\
 &$X_{1-2}^{\epsilon}$ &  &  \\\hline
\textrm{Reissner-Nordstr\"om} &$\left(-\frac{u}{2L^2}-\frac{u^3}{2}\right)$ & $\left(\frac{ku^2}{2L^2}+\frac{2ku^4}{4}\right)$ &
0\\
 &$X_{1-2}^{\epsilon}$ &$\zeta$  & \\\hline
\textrm{Bardeen} &$\left(-\frac{u}{2L^2}-\frac{u^3}{2}\right)$ & 0&$\left(\frac{3ku^5}{4}+\frac{3ku^3}{4L^2}\right)$\\
 &$X_{1-2}^{\epsilon}$ &  & $\zeta$ 
\\\hline\hline
\end{tabular} \end{adjustbox}
\end{table}
\footnotetext{The speed of light is $c=1$, $k$ is a nonzero constant and $\zeta$ denotes that no nontrivial approximate symmetry exists.}
To this end,  the same explicit conditions facilitate the identification of approximate symmetry ``groups" of other critical orbital equations. By way of examples, we list the orbital equations of
the charged rotating Ba\~nados, Teitelboim and Zanelli metric or the  Kerr, Kerr-Newman  or  Kerr-Newman-AdS  spaces, all of which can be rewritten with a perturbation.

\bigskip

\emph{Acknowledgments}

SJ acknowledges financial support from the
National Research Foundation of South Africa (99279). NM acknowledges support from
Standard Bank in partnership with Studytrust. We are grateful to Dr. Andronikos Paliathanasis for his useful comments.


\section{References}

\begin{thebibliography}{99}   

\bibitem{bai2} V.A Baikov, R.K. Gazizov and N.H. Ibragimov,  Approximate symmetries of equations with a small parameter, Mat. Sb. 136 (1988) 435-450 (English Transl. in Math. USSR Sb. 64 (1989) 427-441).

\bibitem{bai3} W.I. Fushchich  and W.M. Shtelen,   On approximate symmetry and approximate solutions of the non-linear wave equation with a small parameter, J. Phys. A: Math. Gen. 22 (1989) 887-890.

      \bibitem{a1} T. Feroze and A.H. Kara, 
  Group theoretic methods for approximate invariants and Lagrangians for some classes of $y''+ \epsilon F(t)y'+y =f(y, y' )$, 
   Int. J. Non-Linear Mech. 37 (2002) 275-280.
      
      \bibitem{a3} A.G. Johnpillai and A.H. Kara, Variational Formulation of Approximate Symmetries and Conservation Laws, 
     Int. J. Theor. Phys. 40 (2001) 1501-1509.     
     
     
  
  
      \bibitem{a2}V.A. Baikov, Approximate symmetries of the van der Pol equation, 
      Differential Equations 30(10) (1994) 1820-1822.

\bibitem{pek1}M. Pakdemirli, M. Y\"ur\"usoy, T. Dolapci, Comparison of Approximate Symmetry Methods for Differential Equations, 
Acta Applic. Math. 80 (2004) 243-271.

\bibitem{bar2} U. Camci, Approximate Noether gauge symmetries of the Bardeen model, Eur. Phys. J. C 74 (2014)  3201. 

\bibitem{gaz} R.K Gazizov, Lie Algebras of Approximate Symmetries, J. Nonlinear Math. Phys., 3 (1996) 96-101.


 \bibitem{r1}  N.H.  Ibragimov, A.H. Kara  and F.M. Mahomed, Lie-Backlund and Noether Symmetries with Applications, 
      Nonl. Dyn. 15 (1998) 115-136.
  \bibitem{r2}  A.G. Johnpillai, A.H. Kara and F.M. Mahomed, Approximate Noether-type symmetries and conservation laws via partial Lagrangians for PDEs with a small parameter,
      J. Comput. Appl. Math. 223 (2009) 508-518.
  \bibitem{r3}  Z.Y. Zhang, Approximate nonlinear self-adjointness and approximate conservation laws, 
      J. Phys. A: Math. Theor. 46  (2013) 155203.
                                                                                      
  \bibitem{r4} A. Paliathanasis and  S. Jamal,  
  Approximate Noether symmetries and collineations for regular perturbative Lagrangians, 
  J. Geom. Phys. {124} (2018)  300-310.

 \bibitem{champ} B. Champagne,  W. Hereman, P. Winternitz, The computer calculation of Lie point symmetries of large systems of differential equations,
 Comp. Phys. Commun.  66 (1991) 319-340.
 
 \bibitem{bau}G. Baumann, Symmetry Analysis of Differential Equations with Mathematica, Springer, New York, 2000. 
  
\bibitem{sym} S. Dimas,  D. Tsoubelis, SYM: A new symmetry-finding package for Mathematica in Group Analysis of Differential Equations, University of Cyprus, Nicosia, Cyprus (2005).

    \bibitem{RNpaper} C.J. Gao, S.N. Zhan,  Reissner-Nordstr\"om metric in the Friedman-Robertson-Walker universe, 
    Phys. Lett. B 595 (2004) 28-35.

\bibitem{g1} M. Maggiore, A generalized uncertainty principle in quantum gravity, Phys. Lett. B 304  (1993) 65-69.
\bibitem{g2}  A. Kempf, Non-pointlike particles in harmonic oscillators, J. Phys. A Math. Gen., 30,  (1997) 2093.
\bibitem{g3} S. Das and E.C. Vagenas, Universality of Quantum Gravity Corrections, Phys. Rev. Lett. 101  (2008) 221301.
\bibitem{g4} S.K. Moayedi, M.R. Setare and H. Moayeri, Quantum Gravitational Corrections to the Real Klein-Gordon Field in the Presence of a Minimal Length, Int. J. Theor. Phys. 49  (2010) 2080.
\bibitem{pal} A. Paliathanasis, S. Pan, S. Pramanik, Scalar field cosmology modified by the Generalized Uncertainty Principle, Class.  Quant. Grav. 32 (24) (2015) 245006. 

\bibitem {StephaniB} H. Stephani, Differential Equations: Their Solutions using
Symmetry, Cambridge University Press, 1989.

\bibitem {Noetherpaper}E. Noether, Nachr. d. K\"{o}nig. Gesellsch. d. Wiss. zu
G\"{o}ttingen, Math-phys. Klasse 235 (1918).


\bibitem{ma}V. Baikov,  R.K. Gazizov, N.H. Ibragimov and F.M. Mahomed,  Closed orbits and their stable symmetries. J. Math.
Phys. 35 (1994) 6525-6535.

\bibitem{o} P.J. Olver,  {Application of Lie Groups to Differential Equations,} Springer, New York, 1993.

\bibitem{prd} N. Dimakis, A. Giacomini,  S. Jamal,  G. Leon and A. Paliathanasis, 
Noether symmetries and stability of ideal gas solutions in
Galileon cosmology,  Phys. Rev. D 95  (2017) 064031.

\bibitem{prd2} S. Capozziello, E. Piedipalumbo, C. Rubano, P. Scudellaro, Noether symmetry approach in phantom quintessence cosmology Phys.Rev. D 80 (2009) 104030.

\bibitem{sit} S. Jamal, A group theoretical application of SO(4,1)  in the de Sitter universe, Gen. Rel. Grav.  49 (88) (2017), DOI 10.1007/s10714-017-2253-4.

\bibitem{olv2} P.J. Olver, { Evolution equations possessing infinitely many symmetries},  
{J. Math. Phys.}  18(6)  (1977)   1212-1215. 

\bibitem{bmm}  S. Jamal, A. Mathebula, Generalized Symmetries and Recursive Operators of Some Diffusive Equations, 
Bull. Malays. Math. Sci. Soc., DOI 10.1007/s40840-017-0510-z.

\bibitem{vit}G.C. McVittie, The mass-particle in an expanding universe, Mon. Not. R. Astron. Soc. 93 (1933) 325-329.

    \bibitem{bar} M. Sharif and S. Waheed,  Energy of Bardeen Model Using Approximate Symmetry Method, 
    Phys. Scr., 83  (2011) 015014.
     \bibitem{sch} A.H. Kara, F. M. Mahomed and A. Qadir,  Approximate symmetries and conservation laws of the geodesic equations for the Schwarzschild metric, 
Nonl. Dyn.  51 (2008) 183-188.
\bibitem{rn} I. Hussain, F. M. Mahomed  and A. Qadir, Second-Order Approximate Symmetries
of the Geodesic Equations for the Reissner-Nordstr\"om Metric and Re-Scaling of Energy of a Test Particle, 
SIGMA  3(115) (2007) 1-9.
      

\end{thebibliography}
\end{document}